\documentclass{article}

\usepackage{amssymb}
\usepackage{mathrsfs}
\usepackage{graphicx}
\usepackage{graphicx}
\usepackage{amsmath}
\usepackage{caption}
\usepackage{flushend,cuted}
\usepackage{color}
\usepackage{float}
\usepackage{subfigure}
\usepackage{multirow}

\newcommand{\be}{\begin{equation}}
\newcommand{\ee}{\end{equation}}
\newcommand{\bit}{\begin{itemize}}
	\newcommand{\eit}{\end{itemize}}

\def\\omega{{\rm Diag}}

\renewcommand{\Re}{{\rm I}\!  {\rm R}}

\newtheorem{proposition}{Proposition}[section]
\newtheorem{theorem}[proposition]{Theorem}

\newtheorem{remark}[proposition]{Remark}

\newtheorem{algorithm}[proposition]{Algorithm}

\newtheorem{proof}[proposition]{proof}

\renewcommand{\Re}{{\rm I}\!  {\rm R}}

\begin{document}
	
	\title{A Majorization Penalty Method for SVM with Sparse {Constraint}}
	
\author{%
	{\sc
		Sitong Lu\thanks{Email: lusitong@bit.edu.cn}} \\[2pt]
	School of Mathematics and Statistics, Beijing Institute of Technology,\\ Beijing, 100081, P. R. China\\[6pt]
	{\sc Qingna Li}\thanks{Corresponding author. Email: qnl@bit.edu.cn. {This author's research is supported by NSF 12071032.}}\\[2pt]
	 School of Mathematics and {Statistics/ Beijing} Key Laboratory on MCAACI,\\ Beijing Institute of Technology,
	Beijing, 100081, P.R.China
}
\maketitle

\begin{abstract}
	{Support vector machine is an important and fundamental technique in machine learning. {Soft-margin SVM models {have} stronger generalization performance compared with {the} hard-margin SVM. Most existing works use the hinge-loss function which can be regarded as an upper bound of the 0-1 loss function. However, it can not explicitly limit the number of misclassified samples. 
	In this paper, we use the idea of soft-margin SVM and propose a new SVM model with {a} sparse constraint. Our model can strictly limit the number of misclassified samples, expressing the soft-margin constraint as a sparse constraint.} {By} constructing a majorization function, a majorization penalty method can be used to solve the {sparse-constrained} optimization problem. We apply {Conjugate-Gradient (CG) method} to solve {the resulting subproblem}. Extensive numerical {results} demonstrate the impressive performance of the proposed majorization penalty method.}
	{Support Vector Machine; Majorization Penalty Method; {Conjugate Gradient Method}; Sparse Constraint.}
\end{abstract}

\section{Introduction}
Support vector machine (SVM) is a traditional and effective machine learning method \cite{awad2015support,vapnik2006estimation,vapnik2013nature,vapnik1996support,yuan2005polynomial,lee2001ssvm}. 
SVM preserves {the} performance benchmark records of handwritten {numerical} recognition, text classification, information retrieval and time series prediction. They are usually used for DNA microarray data analysis \cite{al2006new,labusch2008simple,vapnik1996support,yuan2005polynomial}.
Given {the} training data $x_1,x_2,\dots,x_n\in \Re^m$ and {the corresponding labels} $y_1,y_2,\dots,y_n\in\{-1,1\}$, SVM is aimed to find a hyperplane $\omega^T x + b = 0$ to separate the data according to labels, where $\omega\in \Re^m$ and {the} bias $b\in \Re$ need to be {determined}. The classical {hard-margin} SVM is given as {follows:}
\be\label{classical-svm}
\begin{array}{ll}
	\underset{\omega\in \Re^m, b\in \Re}{\min} & \frac{1}{2}\|\omega\|^2\\
	\ \ \ \ \text{s.t.} 
	& y_i(\omega^T x_i+b)\ge 1 \ \ \ i=1,2,\dots,n.
\end{array}
\ee
In the classical SVM, {we are looking for $(\omega,b)$ such that the hyperplane $\omega^Tx+b=0$ can separate the data as much as possible in the sence that the distince between two types of data is maximized. In such situation,} data are separated strictly, which is based on the assumption
that the data can be linearly separated.

In most cases, optimization methods { are designed to} solve SVM models with soft margin. Some typical examples are {the} \texttt{L1-SVM} and \texttt{L2-SVM} {models}. We will briefly review the methods {for} the above two models, which are related to our work. For a survey of machine learning optimization methods, please refer to \cite{bottou2018optimization,friedman2001elements}.
SVM model with linear summation of {the hinge loss} is very common, which is called \texttt{L1-SVM} as follows
\be\label{l1-svm}
\begin{array}{ll}
	\underset{\omega\in \Re^m, b\in \Re}{\min} & \frac{1}{2}\|\omega\|^2+\rho\sum\limits_{i=1}^{n}\text{max}(1-y_i(\omega^T x_i+b),0).\\
\end{array}
\ee
\texttt{L1-SVM} was put forward by Tibshirani et al. in \cite{zhu20031}. For \texttt{L1-SVM}, {since} {the L1 loss function} is not differentiable, most algorithms cannot be applied. In \cite{fung2004feature}, Mangasarian used Newton's method to solve \texttt{L1-SVM} for selecting features from a very high dimensional space.
In \cite{mangasarian2006exact}, Mangasarian used a generalized Newton's method to solve {the} exact \texttt{L1-SVM}. {Hsieh et al. \cite{hsieh2008dual} proposed a dual coordinate method (DCD) for the dual problem of \texttt{L1-SVM}.} {Yan and Li proposed an augmented Lagrangian method for the primal problem of \texttt{L1-SVM}.
	By using {the} Moreau-Yosida regularization and the proximal operator, they {dealt} with the nonsmooth term of \texttt{L1-SVM} very well \cite{yan2020efficient}.}

If it is the sum of squares {of the hinge loss}, it is called \texttt{L2-SVM}, as shown below:
\be\label{l2-svm}
\begin{array}{ll}
	\underset{\omega\in \Re^m, b\in \Re}{\min} & \frac{1}{2}\|\omega\|^2+\rho\sum\limits_{i=1}^{n}\text{max}(1-y_i(\omega^T x_i+b),0)^2.\\
\end{array}
\ee
The relationship and characteristics of the two models were analyzed in {details} in \cite{burges1999uniqueness,fernandez1998behavior,pontil1998properties,rifkin1999note,yamada2003statistical}.
For the \texttt{L2-SVM} model, Mangasarian \cite{mangasarian2002finite} introduced the finite Newton's method due to the nondifferentiability of the objective function gradient.
It is basically a {unit-step} semismooth {Newton's} method, which uses the inverse of Hessian matrix to calculate {the Newton's} direction. 
{Keerthi and Decoste \cite{keerthi2005modified} proposed an improved Newton's method. They calculated Newton points and performed an accurate line search to determine the step size, which is suited for large scale data mining tasks.
Lin {et al.} porposed the Trust region Newton method (TRN) \cite{lin2007trust} for \texttt{L2-SVM}. Chang et al. \cite{chang2008coordinate} proposed a coordinate descent method for \texttt{L2-SVM}. They {used} Newton's method to solve the subproblem with one variable while fixing other variables. 
{In \cite{hsieh2008dual}, Hsieh et al. also} proposed a dual coordinate descent method (DCD) for the dual {problem} of \texttt{L2-SVM} with large-scale sparse data.
Recently, Hsia et al. \cite{hsia2017study} studied the trust region updating rule in Newton's method. 
Yin and Li {proposed a} semismooth Newton and achieved remarkable results in \cite{yin2019semismooth}.}

\begin{figure}[htbp]
	\centering
	\includegraphics[width=6.7cm,height=5cm]{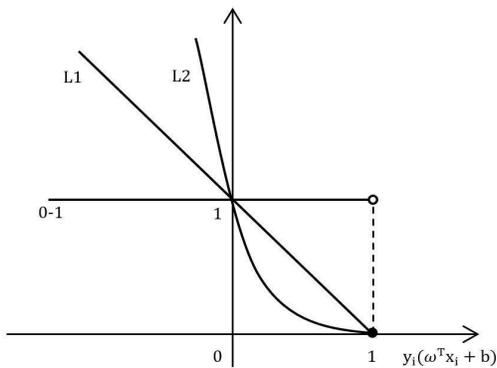}
	\caption{L1, L2 and 0-1 loss function. 
		L1: $\max(0,y_i(\omega^T x_i+b))$; L2: $\max(0,y_i(\omega^T x_i+b))^2$; 0-1: $I(1-y_i(\omega^T x_i+b)>0)$.}
	\label{fig:loss_func}	
\end{figure}

Both \texttt{L1-SVM} and \texttt{L2-SVM} consider the hinge loss function, which is $max(1-y_i(\omega^T x_i+b),0)$. When the sample $(x_i,y_i)$ is correctly classified and the {functional margin} $y_i(\omega^T x_i+b)$ is greater than 1, the loss is 0; otherwise, the loss is $1-y_i(\omega^T x_i+b)$ {or $(1-y_i(\omega^T x_i+b))^2$}. However, the hinge loss {function is} essentially the surrogate of 0-1 loss function, that is {$$I(1-y_i(\omega^T x_i+b)>0):=\left\{
\begin{aligned}
1, & & 1-y_i(\omega^T x_i+b)>0, \\
0, & & \text{otherwise}.\ \ \ \ \ \ \ \ 
\end{aligned}
\right.$$}
We give a figure {drawing the relationship of L1, L2 and 0-1 loss function {in Fig. \ref{fig:loss_func}}.}

Because the 0-1 loss function is discontinuous, few algorithms deal with it directly. The above observations motivate us to propose a new model based on 0-1 loss function. 
The corresponding SVM model with 0-1 loss function is {natually presented as}
\be\label{svm_01loss}
\begin{array}{ll}
	\underset{\omega\in \Re^m, b\in \Re}{\min} &  \frac{1}{2}\|\omega\|^2 +C\sum\limits_{i=1}^n I(1-y_i(\omega^T x_i+b)>0),\\
\end{array}
\ee
where $C>0$ is a regularization parameter. For sample data $i\in\{1,\dots,n\}$, the soft-margin model expects that the sample is correctly classified, i.e., $1-y_i(\omega^T x_i+b)\le0$, {despite a few of them may be misclassified. In other words, {for a very few of} $i\in\{1,2,\dots,n\}$, there is $1-y_i(\omega^T x_i+b)>0$.} Let $s$ $(\gg n)$ be a parameter denoting the {upper bound} {of the} number of samples that can not be classified correctly. Then we can reach the following model 
\be\label{svm_01set}
\begin{array}{ll}
	\underset{\omega\in \Re^m, b\in \Re}{\min} & \frac{1}{2}\|\omega\|^2\\
	\ \ \ \ 	\text{s.t.} & \sum\limits_{i=1}^n I(1-y_i(\omega^T x_i+b)>0) \le s.
\end{array}
\ee
Denote $z=[z_1,z_2,\dots,z_n]^T\in \Re^n$ by $z_i = 1-y_i(\omega^T x_i+b)$. 
Model (\ref{svm_01set}) can be written as the following SVM with {a} sparse constraint
\be\label{svm_sparse}
\begin{array}{ll}
	\underset{\omega\in \Re^m, b\in \Re}{\min} & \frac{1}{2}\|\omega\|^2\\
	\ \ \ \ \text{s.t.} & \| z_+\|_0\le s,\\
	& z_i=1-y_i(\omega^T x_i+b) \ \ \ i=1,2,\dots,n,
\end{array}
\ee
where $(\cdot)_+ :=  \text{max}(\cdot,0)$ and {$\|\cdot\|_0$ is the number of non-zero elements in a vector, {which is often referred to as 0 norm}.} {When $s=0$, (\ref{svm_sparse}) {reduces to} {the} hard-margin SVM problem {(\ref{classical-svm})}.}

0 norm {in (\ref{svm_sparse})} leads to the sparse constrained optimization problem, which is {a} NP-hard problem. In order to deal with the computing challenges of 0 norm, the current algorithm is mainly divided into two mainstreams \cite{zhao2020lagrange}. {One is 'greedy', which means a variety of relaxation forms \cite{zhang2010nearly,gotoh2018dc,zhao2013rsp}. Another is called 'relaxtion', which directly deals with the problem of 0 norm. 
For example, matching pursuit \cite{mallat1993matching} provides a fast and compact method of adaptive function approximation. 
The orthogonal matching pursuit \cite{davis1997adaptive} produces sub-optimal function expansions by iteratively choosing dictionary waveforms.
The subspace pursuit \cite{dai2009subspace} has {{a} low computational complexity and {deals with} very sparse signals.
The compressive sampling matching pursuit (CoSaMP) \cite{needell2009cosamp} offers rigorous bounds
on computational cost and storage.
The hard thresholding pursuit (HTP) \cite{foucart2011hard} combines the iterative hard thresholding (IHT) algorithm
and CoSaMP.
The conjugate gradient iterative hard thresholding \cite{blanchard2014conjugate} combines the low cost of each iteration of the simple line search IHT with the improved convergence speed.}

{{Various} traditional methods have} been extended to sparse constrained nonlinear
optimization. A gradient hard-thresholding method was proposed by Bahmani et al. in \cite{bahmani2013greedy}, which generalizes CoSaMP. Yuan et al. generalized HTP to the sparse constrained convex optimization \cite{yuan2014gradient} and proposed a Newton greedy pursuit method \cite{yuan2014newton}. Zhou et al. \cite{zhou2019global} built the algorithm of Newton Hard-Thresholding Pursuit for the sparse-constrained optimization with the quadratic convergence rate. 
Blumensath and Davies proposed IHT algorithm for the linear compressed sensing problem in \cite{blumensath2008iterative,blumensath2009iterative}. Based on this, Pan et al. \cite{pan2016improved} gave an improved iterative hard thresholding algorithm by employing the Armijo-type stepsize rule, which automatically adjusts the stepsize and {the support} set and leads to a sufficient decrease of the objective function {in} each iteration.

{The inspiration of our paper comes from the framework of majorization projection {method} in \cite{zhou2018fast}}. {For solving the Euclidean distance matrix problem with box constraints {and rank constraint}, it is difficult to {deal with} the matrix rank constraint.} Zhou {et al.} \cite{zhou2018fast} penalized the quadratic distance {function} of {a point to the} conditional positive semidefinite cone with the rank {cut}. By finding the subgradient of {the} penalty function, {they} constructed the majorization function. 
{They also analyzed the convergence of the resulting majorization penalty method.}
Note that the rank of a matrix and the 0 norm of a vector are both {essentially} counting functions. This leads us to think {about} whether we can {apply} the same technique and extend it {to} the vector optimization with {0-norm constraint}. {It is worth noting that there is also work on SVM with 0 norm constraint recently. Zhou established optimality conditions with the stationary equations and solved it by Newton-type method in \cite{zhou2020sparse}.}

{{\bf Our Contribution.} {In this paper, we} first propose a new SVM model, {that is (\ref{svm_01set}),} which considers {0-norm constraint} directly {to restrict the number of misclassified samples no greater than a given level}. We consider the subproblem in penalty method and {penalize} the {0-norm} constraint to the objective function. After that, {a majorization penalty} algorithm framework is proposed, which is also suitable for our SVM model. {The subproblem of the resulting algorithm {is essentially a strongly convex quadratic programming without constraints and} can be efficiently solved by {conjugate gradient} method. Finally, we present extensive numerical results to show the efficiency of the proposed method}.}

This article is organized as follows. {In Section 2, we discus the penalty method for (\ref{svm_sparse})}. In Section 3, we will construct the majorization function,  {and derive the majorization penalty method}. In Section 4, we solve the majorization subproblem by {CG} method. Extensive numerical examples will illustrate the {impressive} performance of our algorithm in Section 5.} {Final conclusions are given in Section 6.}

{{\bf Notation.} 
{Let  
	$ \boldsymbol 1_n\in \Re^n$ be} a vector which all elements are equal to {one}. For a vector $ x \in \Re^n$, $Diag( x)\in \Re^{n\times n}$ is diagonal matrix whose $i$-th diagonal element is $x_i$. We use $\|\cdot\|$ to denote the {L2} norm of vectors.}

\section{Penalty Method for SVM with Sparse Constraint (\ref{svm_sparse})}

In this section, {we consider} the penalty method for SVM with sparse constraint (\ref{svm_sparse}). 

Recall the SVM model with {the} sparse constraint as in (\ref{svm_sparse}). Let $\Omega_s$ be the set defined 
as $\Omega_s:=\{ x\in \Re^n\ |\ \| x_+\|_0\le s\}$. Then (\ref{svm_sparse}) can be equivalently written as follows

\be\label{svm}
\begin{array}{ll}
\underset{\omega\in \Re^m, b\in \Re}{\min} & \frac{1}{2}\|\omega\|^2\\
\ \ \ \ \text{s.t.} & z_i=1-y_i(\omega^T x_i+b), \ \ \ i=1,2,\dots,n,\\
& z\in\Omega_s.
\end{array}
\ee
We call it sparse {constrained} SVM model (SCSVM), where the sparse constraint set $\Omega_s$ is a non-convex set. 

For {the} sparse {constrained set}
$\Omega_s$,
in order to measure the distance between $ z$ and $\Omega_s$, define {$g:\Re^n\rightarrow\Re$} as
$$g( z)=\frac{1}{2}\text{dist}^2( z,\Omega_s),$$
where
\be\label{dist}
\text{dist}( z,\Omega_s):=\min\{\| z- x\|\ |\ x\in\Omega_s\}.\ee 
{Let $\Pi_{\Omega_s}^B(z)$ be the set of {the} optimal {solutions} of {(\ref{dist})}. Due to the nonconvexity of $\Omega_s$, $\Pi_{\Omega_s}^B(z)$ may contain multiple points. Let $\Pi_{\Omega_s}(z)\in\Pi_{\Omega_s}^B(z)$ be one of the elements in $\Pi_{\Omega_s}^B(z)$. Then \be\label{def_g}g( z)=\frac{1}{2}\|z-\Pi_{\Omega_s}(z)\|^2.\ee}
{Therefore, $ z\in\Omega_s$ is equivalent to $g( z)=0$.
{Model} (\ref{svm}) is equivalent to 
\be\label{svm_2}
\begin{array}{ll}
	\underset{\omega\in \Re^m, b\in \Re}{\min} & \frac{1}{2}\|\omega\|^2\\
	\ \ \ \ \text{s.t.} & z_i=1-y_i(\omega^T x_i+b), \ \ \ i=1,2,\dots,n,\\
	& g(z)=0.
\end{array}
\ee
Due to the equivalence between (\ref{svm}) and (\ref{svm_2}), (\ref{svm_2}) is a nonconvex problem and it is in general NP hard. Inspired by \cite{zhou2018fast}, we consider the penalty method to solve (\ref{svm_2}) which penalizes the nonconvex constraint $g(z)=0$.
The subproblem in {the} penalty method is {reformulated as}
\be\label{obj}
\begin{array}{ll}
	\underset{\omega\in \Re^m, b\in \Re}{\min} & F_{\rho}(\omega,b):=\frac{1}{2}\|\omega\|^2+\rho g( z)\\
	\ \ \ \	\text{s.t.} & z_i=1-y_i(\omega^T x_i+ b), \ \ \ i=1,2,\dots,n,
\end{array}
\ee
where $\rho>0$ is a penalty parameter. 

{Denote 
	\be\label{def_PQ}
		A=[ x_1, x_2,\dots, x_n]^T\in \Re^{n\times m}, P=\begin{bmatrix} I_m \ 0 \end{bmatrix} \in \Re^{m\times (m+1)}
		\text{ and }Q=\begin{bmatrix} A \ \boldsymbol 1_n \end{bmatrix}\in \Re^{n\times (m+1)}\ee
	where $n$ is the number of {samples} and $m$ is the number of features.
	{Moreover, let} $\theta\in\Re^{m+1}$ be defined as 
	$$\theta:=\begin{bmatrix} \omega \\ b \end{bmatrix}\in \Re^{m+1}.$$ There is 
	{$$\omega=P\theta \text{ and }
		A\omega+ b=Q\theta.$$}
	The constraints in (\ref{obj})} can be {reformulated as} 
\be\label{def_z} z= \boldsymbol 1_n-Diag(y)(A\omega+ b)= \boldsymbol 1_n-\bar Q\theta,\ee
	where {$\bar Q := Diag(y) Q\in\Re^{n\times(m+1)}$.} 
{Therefore}, (\ref{svm_2}) is equivalent to the following problem
$$
\begin{array}{ll}
\underset{\theta\in \Re^{m+1}}{\min} & f(\theta)\\
\ \ \ \	\text{s.t.} & p(\theta) = 0,
\end{array}
$$
where \be\label{def_f_g} f(\theta):=\frac{1}{2}\|P\theta\|^2 \ \ \text{  and}\ \ \ p(\theta):=g( \boldsymbol 1_n-\bar Q\theta).\ee
{Similarly, substituting (\ref{def_z}) into} $F_{\rho}(\omega,b)$ {in (\ref{obj}), {we get the following} equivalent} new unconstrained problem {of (\ref{obj})}
\be\label{svm-penalty}
\begin{array}{ll}
	\underset{\theta\in \Re^{m+1}}{\min} & F_{\rho}(\theta):=f(\theta)+\rho p(\theta).
\end{array}
\ee
{Now the aim in the penalty method is how we solve the subproblem (\ref{svm-penalty}). To apply majorization method, we will first derive the majorization function for $p(\theta)$ in (\ref{svm-penalty}), as shown in Section 3.}


\section{The Majorizaiton Function}

For subproblem (\ref{svm-penalty}), {it is still not easy to solve. The difficulty lies in $p(\cdot)$. Therefore, similar to as that in \cite{zhou2018fast},}
we would like to construct the majorization function of $p$ {in order to apply the majorization method.} {{To that end}, we introduce the defination of majorization function. The majorization function $p(\theta)$ at {$t\in\Re^{m+1}$}, denoted as $p_m(\theta,t)$, has to satisfy the following conditions
	
	\be\label{majorization} p_m(t,t) = p(t)\ \ \text{and}\ \ p_m(\theta,t)\ge p(\theta)\ \text{for any} \ \theta\in\Re^{m+1}.\ee
	
	We {first} give a theorem {saying} we can have a general choice of majorization function for {the sets which enjoy the orthogonal projection property}.
	
	\begin{theorem}\label{theorem1.1}
		For a {closed} set $\mathcal D$, {let} $\Pi^B_{\mathcal D}( \theta)$ denote the set of projections of a vector $\theta$ {onto} $\mathcal D$ and $\Pi_{\mathcal D}( \theta)\in\Pi^B_{\mathcal D}(\theta)$. Let $k( \theta):=\frac 1 2 \|\Pi_{\mathcal D}( \theta)\|^2$ and $q( \theta):=\frac{1}{2}\| \theta-\Pi_{\mathcal D}( \theta)\|^2$.
		{If} 
		\bit
		{\item[(a)] $\langle\Pi_{\mathcal D}( \theta), \theta-\Pi_{\mathcal D}( \theta)\rangle=0$ for any $\theta$;
			\item[(b)] 
			$\Pi_{\mathcal D}(\theta)\in\partial k(\theta)$, where $\partial k(\theta)$ denotes the subdifferential of $h$ at $\theta$,}
		\eit
		{then for a fixed point $t$, $q_m(\theta,t)$ defined by} $$q_m( \theta, t):=\frac{1}{2}\| \theta\|^2-\langle {\Pi_{\mathcal D}(\theta)}, \theta- t\rangle-k( t)$$ {is a majorization function} of 	
		$q( \theta)$ at $ t$.
	\end{theorem}}
	\begin{proof}
	{By {the} defination of $q(\theta)$, there is}
	\begin{align*}q( \theta)&=\frac{1}{2}\| \theta-\Pi_{\mathcal D}( \theta)\|^2\\
	&=\frac{1}{2}\| \theta\|^2+\frac{1}{2}\|\Pi_{\mathcal D}( \theta)\|^2-\langle\Pi_{\mathcal D}( \theta), \theta\rangle\\
	&{=\frac{1}{2}\| \theta\|^2+\frac{1}{2}\|\Pi_{\mathcal D}( \theta)\|^2-\langle\Pi_{\mathcal D}( \theta), \theta-\Pi_{\mathcal D}( \theta)+\Pi_{\mathcal D}( \theta)\rangle}\\
	&{=\frac{1}{2}\| \theta\|^2-\frac{1}{2}\|\Pi_{\mathcal D}( \theta)\|^2-\langle\Pi_{\mathcal D}( \theta), \theta-\Pi_{\mathcal D}( \theta)\rangle}\\
	&=\frac{1}{2}\| \theta\|^2-{\frac{1}{2}\|\Pi_{\mathcal D}( \theta)\|^2} \ \ \ \ \ \text{(condition (a))}\\
	&\le\frac{1}{2}\| \theta\|^2-\langle {\Pi_{\mathcal D}(\theta)}, \theta- t\rangle-k( t) \ \ \ \ \ \text{(condition (b))}\\
	&=q_m( \theta, t).\end{align*} 
	{Notice {that} $q_m( t, t)=q( t).$} Therefore, {$q_m( \theta, t)$ is a majorization function of $q(\theta)$ at $t$.} \end{proof} 
	
	{In fact,} Theorem \ref{theorem1.1} holds {for} ${\mathcal D}:=\Omega_s$. {To show this, we first give the characterization of $\Pi^B_{\Omega_s}(z).$}
	
	For convenience, for $z\in\Re^n$, we use $z^{\downarrow}\in\Re^n$ to denote the vector whose elements {come from $z$ and} are arranged in a descending order. That is, $z^{\downarrow}_1\ge z^{\downarrow}_2 \ge \dots\ge z^{\downarrow}_n.$ Define 
		\begin{align}\label{alphabeta}\alpha=\{i\ |\ &z_i>z^\downarrow_s\text{ and } z_i>0\}, \ \ \beta=\{i\ |\ z_i=z^\downarrow_s\text{ and } z_i>0\},\\\label{gammatau} &\gamma=\{i\ |\ 0<z_i<z^\downarrow_s\} \text{ and } \tau=\{i\ |\ z_i\le0\}.\end{align}
	
	{Let} $\Pi^B_{\Omega_s}( z)$ denote the set of projections of $z$ {onto} $\Omega_s$.
	Then we give the computation of $\Pi^B_{\Omega_s}(z)$ for $z\in\Re^{n}$ as follows.
	
	\begin{proposition}\label{Pi_def}
		{For any $z\in \Re^n$, there is}
		\be\begin{array}{ll}
			
			\Pi^B_{\Omega_s}(z) = \{x\ |&\ x_i=z_i\ for\ i \in \alpha\cup\beta_1{\cup\tau},\ x_i=0\ for\ i \in (\beta \backslash\beta_1)\cup\gamma,\\
			&where\ \beta_1\subset\beta\ and\ |\beta_1| = s-|\alpha|  \}.
		\end{array}\ee
	\end{proposition}
	\begin{proof}
		Projecting $z$ to {the set} $\Omega_s$ is equivalent to projecting $max(z,0)$ to {the} set $\{x\in\Re^n\ |\ \|x\|_0\le s\}$. {That is, one needs to keep} the largest $s$ components of $max(z,0)$ the same and make the others zero. Note that the negative part of $z$ has no effect on the projection. To sum up, the projection of $z$ to $\Omega_s$ remains components {in} $\alpha\cup\beta_1\cup\tau$, the other components are 0. {Therefore, we get (\ref{Pi_def}). The proof is finished.}\end{proof} 
		
		Due to Proposition \ref{Pi_def}, we have the following results.
	}
	\begin{proposition}\label{corollary1.2}
		Let	$\Pi_{\Omega_s}(z)\in\Pi^B_{\Omega_s}(z)$. {Define $h(\theta):=\frac{1}{2}\|\Pi_{\Omega_s}(\boldsymbol 1_n-\bar Q\theta)\|^2.$} For a {given} vector $\tilde\theta\in\Re^{n+1}$, we have
		\bit
		{\item[(i)] $\langle\Pi_{\Omega_s}( z), z-\Pi_{\Omega_s}( z)\rangle=0$ for any $z{\in\Re^n}$;
			\item[(ii)] 
			$-\bar Q^T\Pi_{\Omega_s}({\boldsymbol 1_n-{\bar Q}\tilde\theta}) \in \partial h(\tilde\theta)$};
		\eit
		
	\end{proposition}
	\begin{proof}
	For (i), we {first show the following holds} \be\label{AAA}(\Pi_{\Omega_s}( z))_i(z-\Pi_{\Omega_s}( z))_i=0,\ i=1,2,\dots,n.\ee 
	If $(\Pi_{\Omega_s}( z))_i=0$, it is obvious {that} $(\Pi_{\Omega_s}( z))_i(z-\Pi_{\Omega_s}( z))_i=0$. If $(\Pi_{\Omega_s}( z))_i \neq 0$, {by Proposition \ref{Pi_def}}, {then $i\in\alpha\cup\beta_1{\cup\tau}$. Therefore,} $(\Pi_{\Omega_s}( z))_i=z_i$. {So $z-\Pi_{\Omega_s}( z)_i=0$ giving again $(\Pi_{\Omega_s}( z))_i(z-\Pi_{\Omega_s}( z))_i=0$. Therefore, (\ref{AAA}) holds.} {To sum up,} {we have} $\langle\Pi_{\Omega_s}( z), z-\Pi_{\Omega_s}( z)\rangle=0$.
	
	For (ii), it is equivalent to the {following} inequality {( letting $\tilde z=\boldsymbol 1_n-\bar Q\tilde\theta$)
			\be\label{subdiff} h(\theta)-h(\tilde\theta)
			\ge \langle-\bar Q^T\Pi_{\Omega_s}(\tilde z),\theta-\tilde\theta\rangle, \forall \  \theta\in\Re^{m+1}.\ee}
	{Note that}
		\begin{align*}
		&\langle-\bar Q^T\Pi_{\Omega_s}(\tilde z),\theta-\tilde\theta\rangle\\
		=&-\langle\Pi_{\Omega_s}(\tilde z),\bar Q(\theta-\tilde\theta)\rangle\\
		=&-\langle\Pi_{\Omega_s}(\tilde z),\bar Q\theta\rangle+\langle\Pi_{\Omega_s}(\tilde z), \bar Q\tilde\theta\rangle\\
		=&-\langle \Pi_{\Omega_s}(\tilde z),\bar Q\theta\rangle
		+\langle \Pi_{\Omega_s}(\tilde z),\boldsymbol 1_n\rangle+
		\langle \Pi_{\Omega_s}(\tilde z), \bar Q\tilde\theta-\boldsymbol 1_n\rangle\\
		=&-\langle \Pi_{\Omega_s}(\tilde z),\bar Q\theta-\boldsymbol 1_n\rangle
		-\langle \Pi_{\Omega_s}(\tilde z),\tilde z \rangle\ \ \ \text{(by $\tilde z=\boldsymbol 1_n-\bar Q\tilde\theta$)}\\
		=&\langle \Pi_{\Omega_s}(\tilde z),z\rangle
		-\|\Pi_{\Omega_s}(\tilde z)\|^2  \ \ \text{(condition (i) and $ z=\boldsymbol 1_n-\bar Q\theta$)}\\
		=&\langle \Pi_{\Omega_s}(\tilde z),z\rangle
		-2h(\tilde\theta).\end{align*}
		So (\ref{subdiff}) is equivalent to {the following conditions:} \be\label{BB}h(\theta)+h(\tilde\theta)\ge \langle \Pi_{\Omega_s}(\tilde z),z\rangle,\ {\forall\ \theta\in\Re^{m+1} \text{ and } z=\boldsymbol 1_n-\bar Q\theta}.\ee
		{To show (\ref{BB}), due} to the {definition} of $\Pi_{\Omega_s}( z)$, we have $$\|\Pi_{\Omega_s}( z)- z\|^2\le\|\Pi_{\Omega_s}(\tilde z)- z\|^2, \ {\forall\ z\in\Re^{n}},$$ {which is equivalent to} 
		$$\|\Pi_{\Omega_s}( z)\|^2-2\langle\Pi_{\Omega_s}( z), z\rangle\le\|\Pi_{\Omega_s}(\tilde z)\|^2-2\langle\Pi_{\Omega_s}(\tilde z), z\rangle.$$
		{By} (i), we have $\langle\Pi_{\Omega_s}( z), z\rangle=\|\Pi_{\Omega_s}( z)\|^2$. {Bringing} it to the above inequality, we get 
		$$-\|\Pi_{\Omega_s}( z)\|^2\le\|\Pi_{\Omega_s}(\tilde z)\|^2-2\langle\Pi_{\Omega_s}(\tilde z), z\rangle.$$
		Rearranging the above in equality gives {(\ref{BB}).}
		{The proof is complete.} \end{proof}
	
	\begin{remark}
		Theorem \ref{theorem1.1} and Proposition \ref{corollary1.2} provide us with a method to construct the majorizaiton funtion of $p(\theta)$ in (\ref{def_f_g}). 
		For fixed $\tilde\theta\in\Re^{m+1}$ and any $\theta\in\Re^{m+1}$, {we have}
		\begin{align*}p(\theta)&=g(z)\\
		&=\frac{1}{2}\| z-\Pi_{\Omega_s}( z)\|^2\\
		&=\frac{1}{2}\| z\|^2+\frac{1}{2}\|\Pi_{\Omega_s}( z)\|^2-\langle\Pi_{\Omega_s}( z), z\rangle\\
		&=\frac{1}{2}\| z\|^2-h(\theta) \ \ \ \ \ \ \text{(Proposition \ref{corollary1.2} (i))}\\
		&\le\frac{1}{2}\| \boldsymbol 1_n-\bar Q\theta\|^2-h(\tilde\theta)+\langle \bar Q^T\Pi_{\Omega_s}({\boldsymbol 1_n-\bar Q\tilde\theta}),\theta- \tilde\theta\rangle\\
		&\ \ \ \ \ \text{ (Proposition \ref{corollary1.2} (ii))}\\
		&=:p_m(\theta, \tilde\theta).\end{align*}
		{Therefore, $p_m(\theta,\tilde\theta)$ is a majorization function of $p(\theta)$ at $\tilde\theta$.}
	\end{remark}
	{
		\section{Majorization Penalty Method} 
		\subsection{Algorithm of Majorization Penalty Method (\texttt{MPM})}
		{Having introduced the majorizaiton function of $p(\theta)$, {we are ready to persent the majorizaiton penalty method. The} idea of majorization penalty method is to solve the penalty problem (\ref{svm-penalty}) by solving the majorization problem. That is, at iteration $\theta^k$, we solve the following  subproblem
			\be\label{svm-majorization}
			\begin{array}{ll}
				\underset{\theta\in \Re^{m+1}}{\min} & F_k(\theta):=\frac{1}{2}\|P\theta\|^2+\rho p_m(\theta,\theta^k).
			\end{array}
			\ee
			The details of majorization penalty method for {solving problem (\ref{svm_sparse})} is given in Algorithm \ref{alg}.}
		
		\begin{algorithm}\label{alg}
			(\texttt{MPM} method {for (\ref{svm_sparse})})
			$$ $$
			\bit
			\item[S1] Input: data matrix $A\in\Re^{n\times m}$, sparse threshold $s>0$, the initial $\theta^0{\in\Re^{m+1}}$ and the penalty parameter $\rho>0$. Set {$k=0$}.
			\item[S2] {Solve {subproblem} (\ref{svm-majorization}) to get a solution {$\theta^{k+1}$}.} 
			\item[S3] If it {reaches the} stopping criterion, {stop; otherwise $k:=k+1$ and go to S2.}
			\eit
		\end{algorithm}
		
		{Next, we mainly discuss the convergence of Algorithm \ref{alg}.} The idea of this proof is basically the same as Themrem 3.7 in \cite{zhou2018fast}.}
	Due to {the} majorization {strategy, there is}
	$$\begin{aligned}
	F_{\rho}(\theta^{k}) &=f(\theta^{k})+\rho p(\theta^{k}) \\
	& =f(\theta^{k})+\rho p_{m}(\theta^{k}, \theta^{k}) \\
	& \geq f(\theta^{k+1})+\rho p_{m}(\theta^{k+1}, \theta^{k}) \\
	& \geq f(\theta^{k+1})+\rho p(\theta^{k+1}) \\
	&=F_{\rho}(\theta^{k+1}).
	\end{aligned}$$
	{As a result,} the seqence $\{F_{\rho}(\theta^k)\}$ generated by Algorithm \ref{alg} is nonincreasing.
	Besides, the {sequence} is bounded below by $0$. Therefore, it converges.
	
	{Next, we will give {two propositions} which will be used later.
	\begin{proposition}\label{nonsingularity} For $\rho>0$,
		$P^TP+\rho Q^TQ$ is a positive definite matrix.
	\end{proposition}
	\begin{proof}
		{For contradiction,} assume that $P^TP+\rho Q^TQ$ is negative semidefinite. There is a nonzero vector $x\in \Re^{m+1}$ such that $$x^T(P^TP+\rho Q^TQ)x\le0.$$ 
		{Note that
		\be\label{PPQQ}P^TP=\begin{bmatrix} I_m & 0\\0 & 0 \end{bmatrix} \ \ \text{and}\ \  Q^TQ=\begin{bmatrix} A^TA & A^T \boldsymbol 1_n\\ \boldsymbol 1_n^TA & n \end{bmatrix}=\begin{bmatrix} \sum\limits_{i=1}^{n}x_ix_i^T & \sum\limits_{i=1}^{n}x_i \\(\sum\limits_{i=1}^{n}x_i)^T & n \end{bmatrix}.\ee
		Therefore}, {we have} $$x^T(P^TP+\rho Q^TQ)x=\sum_{i=1}^{m}x_i^2+\rho x^TQ^TQx\le0.$$ Noticing {that} $Q^TQ$ is a positive semi-definite matrix, the above formula can only be 0. That is $x_i=0$ for $i=1,2,\dots,m$ and $x^TQ^TQx=0$. {On the other hand}, {we have} $$Qx=\begin{bmatrix} A \ \boldsymbol 1_n \end{bmatrix}\begin{bmatrix} 0 \\ \vdots\\0\\ x_{m+1} \end{bmatrix}=x_{m+1}.$$ We get $x_{m+1}=0$. That is, {$x=0$}, which is contradictory. {Therefore, $P^TP+\rho Q^TQ$ is positive definite.} 
	\end{proof}}
{
	
	We need the following assumption.
	
{\bf Assumption 1.} $\tilde{z}\in\Re^n$ has distinct elements. 
	
\begin{proposition}\label{limit} 
	Let $\{ z_k\}$ be a sequence whose limit is $\tilde{ z}$, that is,  $\underset{k\rightarrow\infty}{\lim} z_k=\tilde{ z}$. Suppose $\tilde{z}$ satisfies Assumption 1.	
	Then $\underset{k\rightarrow\infty}{\lim}\Pi_{\Omega_s}( z_k)=\Pi_{\Omega_s}(\tilde{ z})$.
\end{proposition}	
\begin{proof}
	Because we all components of $\tilde{z}$ are distinct, both $\Pi_{\Omega_s}^B( z_k)$ and $\Pi_{\Omega_s}^B(\tilde{ z})$ are singletons. {That is, $\Pi^B_{\Omega_s}( z_k)=\{\Pi_{\Omega_s}( z_k)\}$ and $\Pi^B_{\Omega_s}(\tilde{ z})=\{\Pi_{\Omega_s}(\tilde{ z})\}$.}
	
	Donate $\alpha(z)$ and $\beta(z)$ as {in} (\ref{alphabeta}) with {respect to} $z$.
	First, we will prove that $$\alpha(z_k)\cup\beta(z_k)=\alpha(\tilde{ z})\cup\beta(\tilde{ z})$$ when $k$ is large enough. 
	Suppose the above equation does not hold, which means that there is an index $i$ such that $i\in\alpha(z_k)\cup\beta(z_k)$ and $i\notin\alpha(\tilde{ z})\cup\beta(\tilde{ z})$ or an index $j$ such that $j\notin\alpha(z_k)\cup\beta(z_k)$ and $j\in\alpha(\tilde{ z})\cup\beta(\tilde{ z})$. Because $|\alpha(z_k)\cup\beta(z_k)|=|\alpha(\tilde{ z})\cup\beta(\tilde{ z})|=s$\footnote{Here we use $|\cdot|$ to donate the number of elements in a set.}, such $i$ and $j$ exist at the same time. Above all, we get two sets of inequalities as follows
	$$(z_k)_i\ge(z_k)_s^{\downarrow}>(z_k)_j \ \text{ and }\ (\tilde{ z})_i<(\tilde{ z})_s^{\downarrow}\le(\tilde{ z})_j.$$ Combined {the  above formulae} with $\underset{k\rightarrow\infty}{\lim} z_k=\tilde{ z}$ and the order preserving property of sequence limit, it is contradictory. So $\alpha(z_k)\cup\beta(z_k)=\alpha(\tilde{ z})\cup\beta(\tilde{ z})$. {Therefore,} we get $\underset{k\rightarrow\infty}{\lim}(\Pi_{\Omega_s}( z_k))_i=(\Pi_{\Omega_s}( \tilde{ z}))_i$ for any $i\in\alpha(z_k)\cup\beta(z_k)=\alpha(\tilde{ z})\cup\beta(\tilde{ z})$.
	
	For $i\notin\alpha(\tilde{ z})\cup\beta(\tilde{ z})$,  note that
	$(\Pi_{\Omega_s}( z_k))_i=min\{( z_k)_i,0\}$ and $(\Pi_{\Omega_s}( \tilde{ z}))_i=min\{\tilde{ z}_i,0\}$. Because $min\{\cdot,0\}$ is a continuous monotone increasing function and $\underset{k\rightarrow\infty}{\lim} z_k=\tilde{ z}$, we get $\underset{k\rightarrow\infty}{\lim}(\Pi_{\Omega_s}( z_k))_i=(\Pi_{\Omega_s}( \tilde{ z}))_i$ for any $i\notin\alpha(z_k)\cup\beta(z_k)$.
	
	Therefore, $\underset{k\rightarrow\infty}{\lim}\Pi_{\Omega_s}( z_k)=\Pi_{\Omega_s}(\tilde{ z})$ holds. The proof is finished.
\end{proof}}

{
	Based on Proposition \ref{nonsingularity} and Proposition \ref{limit}}, we will {a} give stronger convergence {below}.
	\begin{theorem}\label{convergence}
		{Let} $\{\theta^k\}$ be the sequence generated by the \texttt{MPM}. We have
{the following results.}
\begin{itemize}		
		\item[(i)] $F_{\rho}(\theta^{k+1})-F_{\rho}(\theta^{k})\le-\frac{\rho}{2}\|\theta^{k+1}-\theta^{k}\|^2,\ k=1,2\dots$.
		Conseqencetly, $\|\theta^{k+1}-\theta^{k}\|\rightarrow0.$
		
		\item[(ii)] Let $\hat{\theta}$ be an accumulation point of $\{\theta^k\}$ {and assume that $\hat{z}=\boldsymbol 1_n-\bar Q\hat\theta$ satisfy Assumption 1.} Then for any $\theta$, we have	
		{$$\langle\nabla f(\hat{\theta})-\rho\bar Q^T {(\boldsymbol 1_n-\bar Q\hat{\theta})}+\rho\bar Q\Pi_{\Omega_s}(\boldsymbol 1_n-\bar Q\hat{\theta}), \theta-\hat{\theta}\rangle\ge0,$$} whcih means that $\hat{\theta}$ is a stationary point of (\ref{svm-penalty}).
		
		{\item[(iii)] If $\hat{\theta}$ is an isolated accumulation point of the sequence $\{\theta^k\}$, then the whole sequence $\{\theta^k\}$ converges to $\hat{\theta}$.}
\end{itemize}
	\end{theorem}
\begin{proof}
(i) {Let}
\be\label{CC}z^k:=\boldsymbol 1_n-{\bar Q}\theta^k \text{ and } z^{k+1}:=\boldsymbol 1_n-{\bar Q}\theta^{k+1}.\ee
{Note that} 
{\be\label{f-f} f(\theta^{k+1})-f(\theta^k)=\frac{1}{2}\langle P(\theta^{k+1}+\theta^{k}),P(\theta^{k+1}-\theta^{k})\rangle=\frac{1}{2}\langle P^TP(\theta^{k+1}+\theta^{k}),\theta^{k+1}-\theta^{k}\rangle\ee}
{and}
{\be\label{nabla_f}\nabla f(\theta^{k+1})=P^TP\theta^{k+1}.\ee }
{Moreover, there is}
{$$\begin{aligned}\label{convergence_2} \|z^{k+1}\|^{2}-\|z^{k}\|^{2} =2\langle z^{k+1}-z^{k}, z^{k+1}\rangle-\|z^{k+1}-z^{k}\|^{2}\\
		=2\langle \theta^{k+1}-\theta^{k}, -\bar Q^{T} z^{k+1}\rangle-\|z^{k+1}-z^{k}\|^{2}.\end{aligned}$$}
{By} {Proposition} \ref{corollary1.2} (ii), {we have}
\be\label{convergence_3} h(\theta^{k+1})-h(\theta^{k})\ge\langle -\bar Q^T\Pi_{\Omega_s}( {z^k}),\theta^{k+1}-\theta^{k}\rangle.\ee
{Since $\theta^{k+1}$ is the optimal solution of (\ref{svm-majorization}), by} the optimality condition of (\ref{svm-majorization}), there is \be\label{convergence_4}\nabla F_k(\theta^{k+1})=0,\ee
	where $\nabla F_k(\theta^{k+1})=\nabla f(\theta^{k+1})-\rho \bar Q^T(\boldsymbol 1_n-{\bar Q}\theta^{k+1})+\rho\bar Q^T\Pi_{\Omega_s}(z^k)$.

Therefore, {$$\begin{aligned}
	& F_{\rho}(\theta^{k+1})-F_{\rho}(\theta^{k}) \\
	=& f(\theta^{k+1})-f(\theta^{k})+\rho p(\theta^{k+1})-\rho p(\theta^{k}) \ \ (\text{by} (\ref{def_f_g}))\\
	=& \frac{1}{2}\langle P^TP(\theta^{k+1}+\theta^k),\theta^{k+1}-\theta^k\rangle
	+(\rho / 2)(\left\|z^{k+1}\right\|^{2}-\left\|z^{k}\right\|^{2})\\
	&-\rho(h(\theta^{k+1})-h(\theta^{k})) \ \ (\text{by } (\ref{f-f}))\\
	=& \langle\nabla f(\theta^{k+1}),\theta^{k+1}-\theta^k\rangle-\frac{1}{2}\langle P^TP(\theta^{k+1}-\theta^k),\theta^{k+1}-\theta^k\rangle\\
	+&(\rho / 2)(\left\|z^{k+1}\right\|^{2}-\left\|z^{k}\right\|^{2})-\rho(h(\theta^{k+1})-h(\theta^{k})) \ \ (\text{by } (\ref{nabla_f}))\\
	=&\left\langle\nabla f(\theta^{k+1})-\rho \bar Q z^{k+1}, \theta^{k+1}-\theta^{k}\right\rangle-\frac 1 2 \|P(\theta^{k+1}-\theta^k)\|^2 \\
	&-(\rho / 2)\left\|z^{k+1}-z^{k}\right\|^{2}-\rho(h(\theta^{k+1})-h(\theta^{k})) \ \ (\text{by } (\ref{convergence_2}))\\
	\leq &\left\langle\nabla f(\theta^{k+1})-\rho \bar Q z^{k+1}+\rho \bar Q^T\Pi_{\Omega_s}(z^{k}), \theta^{k+1}-\theta^{k}\right\rangle \\
	&-\frac 1 2\|P(\theta^{k+1}-\theta^k)\|^2-(\rho / 2)\left\|z^{k+1}-z^{k}\right\|^{2}\ \ (\text{by } (\ref{convergence_3}))\\
	=& -\frac 1 2\|P(\theta^{k+1}-\theta^k)\|^2 -(\rho / 2)\left\|\bar Q(\theta^{k+1}-\theta^{k})\right\|^{2}\ \ (\text{by } (\ref{convergence_4}))\\
	=& -\frac 1 2\langle(P^TP+\rho Q^TQ)(\theta^{k+1}-\theta^{k}),\theta^{k+1}-\theta^{k}\rangle.
	\end{aligned}$$}
{Due to Proposition \ref{nonsingularity}, $P^TP+\rho Q^TQ$ is a positive definite matrix.
{Let $\lambda_{min}>0$ be the smallest eigenvalue of $P^TP+\rho Q^TQ.$ There is}	
$$F_{\rho}(\theta^{k+1})-F_{\rho}(\theta^{k})\le-\frac 1 2\lambda_{min}\|\theta^{k+1}-\theta^k\|^2.$$}
This proves that the sequence $\{F_{\rho}(\theta^k)\}$ is non-increasing and it is also bounded below by 0.
Taking the limits on both sides, we have $\|\theta^{k+1}-\theta^{k}\|^2\rightarrow 0$.

(ii) We first prove that $\{\theta_k\}$ generated by Algorithm \ref{alg} is bounded. Because $\{F_{\rho}(\theta^k)\}$ is non-increasing and it is also bounded below by 0, {as} proved in (i), $\{F_{\rho}(\theta^k)\}$ is bounded. {Therefore, we have} $$F_{\rho}(\theta^k)=\frac 1 2 \|\omega_k\|^2+\frac\rho 2\|z_k-\Pi_{\Omega_s}(z_k)\|^2\ge\frac 1 2 \|\omega_k\|^2,$$ {which} shows that both $\{\omega_k\}$ and $\{\|z_k-\Pi_{\Omega_s}(z_k)\|\}$ are bounded. The latter means that {for any index $i\in(\beta(z_k)\ \backslash\beta_1(z_k))\cup\gamma(z_k)$, there is $\{(z_k)_i\}$ is} bounded, which leads to the boundness of $\{(z_k)_i\}$. {Note that $(z_k)_i=1-y_i(\omega_k^Tx_i+b_k)$, which means that $b_k$ is bounded.} So we proved that $\{\theta_k\}$ is bounded.

Therefore, we suppose $\hat{\theta}$ is the limit of a subsequence{, denoted by} $\{\theta^{k_l}\}$. Due to $\|\theta^{k+1}-\theta^k\|^2\rightarrow 0$ in (i), we get the sequence $\{\theta^{k_l+1}\}$ also converges to $\hat{\theta}$.
 Since $\theta^{k+1}$ is the optimal solution of (\ref{svm-majorization}), by the optimality condition of (\ref{svm-majorization}), there is $$\nabla F_k(\theta^{k+1})^T(\theta-\theta^{k+1})\ge0,$$
 {That is,}
\be\label{convergence_5} \langle\nabla f(\theta^{k+1})-\rho Q^T z^{k+1}+\rho Q^T\Pi_{\Omega_s}(z^{k}), \theta-\theta^{k+1}\rangle \geq 0.\ee

{Due to Proposition \ref{limit}, we can take} the limit for both sides of (\ref{convergence_5}) {and get $$\langle\nabla f(\hat{\theta})-\rho\bar Q^T (\boldsymbol 1_n-\bar Q\hat{\theta})+\rho\bar Q\Pi_{\Omega_s}(\boldsymbol 1_n-\bar Q\hat{\theta}), \theta-\hat{\theta}\rangle\ge0.$$ Therefore, $\hat{\theta}$ is a stationary point of (\ref{svm-penalty}).} 

{(iii) We have proved $\|\theta^{k+1}-\theta^{k}\|\rightarrow0$ in (i). From Theorem 5.4 in \cite{kanzow1999qp}, the whole sequence converges to $\hat{\theta}$.}
\end{proof}

\subsection{Solving Subproblem (\ref{svm-majorization})}

In (\ref{svm-majorization}), {the function $F_k(\theta)$ can be simplified as} \begin{align*}F_k(\theta)&=\frac{1}{2}\|P\theta\|^2+\frac{\rho}{2}\| \boldsymbol 1_n-\bar Q\theta\|^2+\rho\langle\Pi_{\Omega_s}(\boldsymbol 1_n-\bar Q\theta^k), \bar Q(\theta-\theta^k)\rangle-\rho h(\theta^k)\\
&=\frac{1}{2}\|P\theta\|^2+\frac{\rho}{2}\| \boldsymbol 1_n-\bar Q\theta\|^2+\rho\langle\Pi_{\Omega_s}(\boldsymbol 1_n-\bar Q\theta^k), \bar Q\theta\rangle+C_0\\
&=\frac{1}{2}\|P\theta\|^2+\frac{\rho}{2}\|\bar Q\theta\|^2+\rho\langle\Pi_{\Omega_s}(\boldsymbol 1_n-\bar Q\theta^k)- \boldsymbol 1_n,\bar Q\theta\rangle+C_1\\
&=\frac{1}{2}\|P\theta\|^2+\frac{\rho}{2}\|\bar Q\theta\|^2+\rho\langle \bar Q^T\Pi_{\Omega_s}(\boldsymbol 1_n-\bar Q\theta^k)- \boldsymbol 1_n,\theta\rangle+C_1,
\end{align*}
{ where $C_0=\rho\langle\Pi_{\Omega_s}(\boldsymbol 1_n-\bar Q\theta^k), \bar Q\theta^k\rangle-\rho h(\theta^k$) and $C_1=C_0+\frac \rho 2\|\boldsymbol 1_n\|^2$ denote the constant {parts} with respect to $\theta$.} 

{Notice that $F_k(\theta)$ is a convex function, solving (\ref{svm-majorization}) is equivalent to solve \be\label{nabla}\nabla F_k(\theta)=0,\ee} {where}
$$
\nabla F_k(\theta)=(P^TP+\rho Q^TQ)\theta+\rho \bar Q^T(\Pi_{\Omega_s}(\boldsymbol 1_n-\bar Q\theta^k)- \boldsymbol 1_n).
$$
{Equation (\ref{nabla}) can be regarded as}
\be\label{Newton}(P^TP+\rho Q^TQ)\theta=\rho \bar Q^T(\boldsymbol 1_n-\Pi_{\Omega_s}(\boldsymbol 1_n-\bar Q\theta^k)),\ee
{where $P^TP$ and $Q^TQ$ are defiend as in (\ref{PPQQ}).}


{Note that by} {Proposition \ref{nonsingularity},} $P^TP+\rho Q^TQ$ is a positive definite matrix, which {actually shows that (\ref{svm-majorization}) is an unconstrained strongly convex quadratic programming problem. It} ensures {that solving (\ref{Newton}) is meaningful} {and one will get an unique optimal solution of (\ref{svm-majorization}).}

To sum up, {in majorization penalty method \texttt{MPM}, solving each subproblem (\ref{svm-majorization}) reduces to solving a linear system (\ref{nabla}), which has a symmetric positive definite coefficient matrix. To further speed up the algorithm, below we discuss more details about solving (\ref{nabla}) more efficiently.}
{
\begin{remark}
	For {(\ref{nabla})},
	when $m$ is small ($m<100$), we can {compute} the inverse of a matrix directly. {That means, we solve the linear system (\ref{svm-majorization}) by direct solver.} {In other cases}, {We use conjugate gradient method, which is {an iterative solver and is suitable for {solving} large-scale linear systems}.}
\end{remark}}

\section{Numberical Experiment and Comparsion}
In this section, we first discuss the performance of our algorithm with different choices of parameters. Then we will compare it with the {state-of-art} algorithms provided by {LIBLINEAR in} \cite{fan2008liblinear}.

All experiments are tested in Matlab R2016b in Windows 10 on a Lenovo
desktop computer with an Intel(R) Core(TM) i5-9300H CPU at 2.40 GHZ
and 16 GB of RAM. 

For the stopping criterion, {similar to that in \cite{zhou2018fast}, we define two residuals} to measure the convergence of $f$ and $p$ respectively. {That is,}
$$\texttt{f-prog}_k:=\frac{f(\theta^{k-1})-f(\theta^{k})}{\rho+f(\theta^{k-1})}$$
and
$$\begin{aligned}\texttt{p-prog}_k:&=\frac{2p(\theta^{k})}{\|\theta^k\|^2}=\frac{2g( \boldsymbol 1_n-\bar Q\theta^{k})}{\|\theta^{k}\|^2}.\end{aligned}$$
The smaller $\texttt{p-prog}_k$ means {that} $z$ is closer to $\Omega_s$ in (\ref{svm}). 
\texttt{MPM} stops when {the following conditions hold} $$\texttt{f-prog}_k\le\sqrt{n}10^{-3} \ \text{and} \ \texttt{p-prog}_k\le10^{-3}.$$ 

\subsection{Performance of Our Algorithm}
In this {part},
we use \texttt{LIBLINEAR}'s standard real data set for classification ({$38$} data sets). For some classification {data sets} whose labels do not belong to $\{-1,1\}$, we change their labels and set them to belong to $\{-1,1\}$. For example, for the dataset {"breast-cancer"}, the label of the sample is $2$ or $4$. {We} turn tag $2$ into $1$ and tag $4$ into $1$. Similarly, we use the same strategy for the dataset: {"liver-disease"}, "mushrooms", "phishing" and "svmguide1". The detailed information of classification data sets is shown in Table \ref{tab:dataset}.

\begin{table}[tbhp]
{\footnotesize
	\caption{Data Information for Classification ($n$ is the number of instances, $m$ is
		the number of features, $\#$nonzeros represents the number of non-zero elements
		in all training instances and density shows the ratio: $\#$nonzeros/(m·n)).} \label{tab:dataset}
	\setlength{\tabcolsep}{5mm}{
		\begin{tabular}{ccccc} \hline
			Data set & n & m & $\#$nonzeros & density \\ \hline
			a1a & 30956 & 123 & 429343 & 11.28$\%$ \\
			a2a & 30296 & 123 & 420188 & 11.28$\%$ \\
			a3a & 29376 & 123 & 407430 & 11.28$\%$ \\
			a4a & 27780 & 123 & 385302 & 11.28$\%$ \\
			a5a & 26147 & 123 & 362653 & 11.28$\%$ \\
			a6a & 21341 & 123 & 295984 & 11.28$\%$ \\
			a7a & 16461 & 123 & 228288 & 11.28$\%$ \\
			a8a & 22696 & 123 & 314815 & 11.28$\%$ \\
			a9a & 32561 & 123 & 451592 & 11.28$\%$ \\
			w1a & 47272 & 300 & 551176 & 3.89$\%$ \\
			w2a & 46279 & 300 & 539213 & 3.89$\%$ \\
			w3a & 44837 & 300 & 522338 & 3.89$\%$ \\
			w4a & 42383 & 300 & 493583 & 3.89$\%$ \\
			w5a & 39861 & 300 & 464466 & 3.89$\%$ \\
			w6a & 32561 & 300 & 379116 & 3.89$\%$ \\
			w7a & 25057 & 300 & 291438 & 3.89$\%$ \\
			w8a & 49749 & 300 & 579586 & 3.89$\%$ \\
			australian & 690 & 14 & 8447 & 87.44$\%$ \\
			breast-cancer & 638 & 10 & 6380 &100$\%$ \\
			cod-rna & 59535 & 8 & 476280 & 100$\%$ \\
			colon-cancer & 62 & 2000 & 124000 & 100$\%$ \\
			diabetes & 768 & 8 & 6135 & 99.85$\%$ \\
			duke breast-cancer & 38 & 7129 & 270902 & 100$\%$ \\
			fourclass & 862 & 2 & 1717 & 99.59$\%$ \\
			german.numer & 1000 & 24 & 23001 & 95.84$\%$ \\
			heart & 270 & 13 & 3510 & 100$\%$ \\
			ijcnn1 & 49990 & 22 & 649870 & 59.09$\%$ \\
			ionosphere & 351 & 34 & 10551 & 88.41$\%$ \\
			leukemia & 38 & 7129 & 270902 & 100$\%$ \\
			liver-disorders & 145 & 5 & 725 & 100$\%$ \\
			mushrooms & 8124 & 112 & 170604 & 18.75$\%$ \\
			phishing & 11055 & 68 & 331610 & 44.11$\%$ \\
			skin nonskin & 245057 & 3 & 735171 & 100$\%$ \\
			splice & 2175 & 60 & 130500 & 100$\%$ \\
			sonar & 208 & 60 & 12479 & 99.99$\%$ \\
			svmguide1 & 3089 & 4 & 12356 & 100$\%$ \\
			svmguide3 & 1243 & 22 & 27208 & 99.50$\%$ \\
			covtype.binary & 581012 & 54 & 6940438 & 22.12$\%$ \\ \hline
			
		\end{tabular}}
	}
\end{table}

We {choose the parameters as follows:} the initial penalty parameter $\rho=0.4$ in Algorithm 1 Step 3, the initial point $\theta^0=0$, maximum iteration {number} of CG method is $500$. {The CG method is stopped if the residual in CG is smaller than $10^{-3}$.}  We report the number of iterations $k$, the total number of CG
iterations $cg$, the cputime $t$ and the accuracy of the algorithm with different sparse ratios $SR$, where 
the accuracy is calculated by 
$$\frac{\text{number of test data whose predicted labels are correct}}{\text{number of test {data}}}$$
and {the sparse ratio $SR$ can be calculated with 
	$$SR:=\frac s n \times 100\%.$$ }
In Table \ref{tab:test_algo} {and Table \ref{tab:test_algo2}, we use $'0'$ to indicate} that we did not use CG method in this case. 

From Table \ref{tab:test_algo} {and Table \ref{tab:test_algo2}}, we get the following conclusions:
{(i)} The algorithm has solved all the problems successfully. This shows that our method can solve the problem (\ref{svm}) in a short cuptime.
{(ii)} When the $SR$ is very high (such as $50\%$), the cputime $t$, the total number of CG
iterations $cg$ and the number of majorization iterations $k$ {are} lower than other cases. This is because when the $SR$ is very {high}, it can be regarded as {we} no longer considering the sparse constraint, and the algorithm can converge faster. {In other words, when the $SR$ is very high}, the convergence criterion is easier to achieve.
(iii) From the result of accuracy, {the same} data set has different solutions with different $SR$s. It is difficult to find a {suitable $SR$} for {all} data sets. In the following numerical experiments, {we set $SR$ to $10\%$.}

\begin{table}[tbhp]
{\footnotesize
	\caption{ {Numerical Results for \texttt{MPM} with Different Sparse Ratio (I).}} \label{tab:test_algo} 
	\begin{center}
			\begin{tabular}{p{3mm}p{2.5mm}p{2mm}p{2.5mm}p{5mm}p{12mm}|p{12mm}p{2.5mm}p{2mm}p{2.5mm}p{5mm}p{12mm}} \hline
				data set &$SR$& k & $cg$ & t(s) & accuracy & data set & $SR$ & k & $cg$ & t(s) & accuracy \\ \hline
				\multirow{5}{*}{a1a} & 1$\%$ & 17 & 1061 & 0.545 & 84.1775$\%$ & \multirow{5}{*}{a2a} & 1$\%$ & 16 & 1152 & 0.027 & 84.3346$\%$\\
				~ & 5$\%$ & 24 & 1486 & 0.035 & 83.7253$\%$ & ~ & 5$\%$ & 22 & 1596 & 0.035 & 84.1827$\%$\\
				~ & 10$\%$ & 27 & 1699 & 0.037 & 83.4992$\%$ & ~ & 10$\%$ & 25 & 1811 & 0.040 & 84.1926$\%$\\
				~ & 15$\%$ & 23 & 1409 & 0.034 & 83.5024$\%$ & ~ & 15$\%$ & 23 & 1668 & 0.038 & 84.0045$\%$\\
				~ & 25$\%$ & 21 & 1303 & 0.032 & 80.4917$\%$ & ~ & 25$\%$ & 16 & 1203 & 0.025 & 81.4266$\%$\\
				~ & 50$\%$ & 11 & 729 & 0.017 & 83.2730$\%$ & ~ & 50$\%$ & 11 & 836 & 0.018 & 83.4599$\%$\\\hline	
				\multirow{5}{*}{a3a} & 1$\%$ & 18 & 1487 & 0.035 & 84.3239$\%$ & \multirow{5}{*}{a4a} & 1$\%$ & 19 & 1791 & 0.052 & 84.4888$\%$\\
				~ & 5$\%$ & 25 & 2056 & 0.049 & 84.1435$\%$ & ~ & 5$\%$ & 25 & 2420 & 0.061 & 84.4240$\%$\\
				~ & 10$\%$ & 25 & 2083 & 0.043 & 83.7078$\%$ & ~ & 10$\%$ & 26 & 2574 & 0.059 & 84.3269$\%$\\
				~ & 15$\%$ & 24 & 1970 & 0.041 & 83.2755$\%$ & ~ & 15$\%$ & 26 & 2603 & 0.060 & 83.5421$\%$\\
				~ & 25$\%$ & 17 & 1437 & 0.030 & 80.4739$\%$ & ~ & 25$\%$ & 18 & 1858 & 0.048 & 81.3139$\%$\\
				~ & 50$\%$ & 11 & 976 & 0.023 & 83.2006$\%$ & ~ & 50$\%$ & 12 & 1255 & 0.034 & 83.8013$\%$\\\hline	
				\multirow{5}{*}{a5a} & 1$\%$ & 20 & 2129 & 0.054 & 84.5642$\%$ & \multirow{5}{*}{a6a} & 1$\%$ & 19 & 2451 & 0.078 & 84.7383$\%$\\
				~ & 5$\%$ & 25 & 2701 & 0.071 & 84.5298$\%$ & ~ & 5$\%$ & 25 & 3208 & 0.097 & 84.7196$\%$\\
				~ & 10$\%$ & 27 & 2912 & 0.067 & 84.3653$\%$ & ~ & 10$\%$ & 26 & 3358 & 0.113 & 84.5368$\%$\\
				~ & 15$\%$ & 25 & 2732 & 0.064 & 83.7763$\%$ & ~ & 15$\%$ & 24 & 3085 & 0.097 & 83.3372$\%$\\
				~ & 25$\%$ & 18 & 2034 & 0.053 & 80.8965$\%$ & ~ & 25$\%$ & 20 & 2632 & 0.088 & 79.4574$\%$\\
				~ & 50$\%$ & 12 & 1297 & 0.035 & 83.4551$\%$ & ~ & 50$\%$ & 12 & 1576 & 0.062 & 83.6652$\%$\\\hline	
				\multirow{5}{*}{a7a} & 1$\%$ & 19 & 2712 & 0.097 & 84.9098$\%$ & \multirow{5}{*}{a8a} & 1$\%$ & 19 & 2850 & 0.120 & 85.3421$\%$\\
				~ & 5$\%$ & 25 & 3527 & 0.124 & 84.9341$\%$ & ~ & 5$\%$ & 26 & 3990 & 0.149 & 85.1698$\%$\\
				~ & 10$\%$ & 26 & 3678 & 0.131 & 84.7640$\%$ & ~ & 10$\%$ & 27 & 4127 & 0.157 & 85.2712$\%$\\
				~ & 15$\%$ & 25 & 3575 & 0.121 & 83.7191$\%$ & ~ & 15$\%$ & 29 & 4485 & 0.171 & 83.7304$\%$\\
				~ & 25$\%$ & 16 & 2361 & 0.085 & 80.9611$\%$ & ~ & 25$\%$ & 18 & 2790 & 0.109 & 79.9088$\%$\\
				~ & 50$\%$ & 11 & 1586 & 0.066 & 83.4639$\%$ & ~ & 50$\%$ & 11 & 1658 & 0.076 & 84.1257$\%$\\\hline	
				\multirow{5}{*}{a9a} & 1$\%$ & 19 & 3269 & 0.153 & 85.0378$\%$ & \multirow{5}{*}{w1a} & 1$\%$ & 25 & 1177 & 0.045 & 96.8668$\%$\\
				~ & 5$\%$ & 25 & 4247 & 0.197 & 84.8658$\%$ & ~ & 5$\%$ & 8 & 436 & 0.019 & 95.1873$\%$\\
				~ & 10$\%$ & 27 & 4591 & 0.196 & 84.7307$\%$ & ~ & 10$\%$ & 10 & 524 & 0.022 & 95.1736$\%$\\
				~ & 15$\%$ & 25 & 4252 & 0.185 & 83.5821$\%$ & ~ & 15$\%$ & 10 & 525 & 0.022 & 95.1668$\%$\\
				~ & 25$\%$ & 17 & 2911 & 0.141 & 78.9755$\%$ & ~ & 25$\%$ & 9 & 483 & 0.021 & 93.7370$\%$\\
				~ & 50$\%$ & 12 & 2005 & 0.117 & 83.2320$\%$ & ~ & 50$\%$ & 14 & 711 & 0.030 & 93.2718$\%$\\\hline	
				\multirow{5}{*}{w2a} & 1$\%$ & 32 & 4404 & 0.282 & 98.2889$\%$ & \multirow{5}{*}{w3a} & 1$\%$ & 31 & 1018 & 0.036 & 97.7519$\%$\\
				~ & 5$\%$ & 9 & 1530 & 0.124 & 96.8660$\%$ & ~ & 5$\%$ & 14 & 551 & 0.020 & 97.6894$\%$\\
				~ & 10$\%$ & 7 & 1200 & 0.113 & 95.9926$\%$ & ~ & 10$\%$ & 14 & 586 & 0.021 & 97.5110$\%$\\
				~ & 15$\%$ & 9 & 1473 & 0.120 & 96.0847$\%$ & ~ & 15$\%$ & 9 & 377 & 0.015 & 97.0738$\%$\\
				~ & 25$\%$ & 9 & 1530 & 0.127 & 95.4965$\%$ & ~ & 25$\%$ & 8 & 318 & 0.013 & 96.9200$\%$\\
				~ & 50$\%$ & 9 & 1572 & 0.130 & 92.9418$\%$ & ~ & 50$\%$ & 10 & 392 & 0.016 & 96.4784$\%$\\\hline	
				\multirow{5}{*}{w4a} & 1$\%$ & 28 & 2051 & 0.096 & 97.1965$\%$ & \multirow{5}{*}{w5a} & 1$\%$ & 23 & 2096 & 0.105 & 98.4697$\%$\\
				~ & 5$\%$ & 7 & 654 & 0.037 & 94.7947$\%$ & ~ & 5$\%$ & 8 & 879 & 0.057 & 96.9820$\%$\\
				~ & 10$\%$ & 8 & 719 & 0.040 & 94.7660$\%$ & ~ & 10$\%$ & 9 & 949 & 0.058 & 96.9695$\%$\\
				~ & 15$\%$ & 10 & 869 & 0.047 & 94.6635$\%$ & ~ & 15$\%$ & 10 & 1048 & 0.066 & 96.9619$\%$\\
				~ & 25$\%$ & 8 & 722 & 0.048 & 92.6756$\%$ & ~ & 25$\%$ & 8 & 867 & 0.055 & 95.7026$\%$\\
				~ & 50$\%$ & 12 & 1006 & 0.052 & 91.6510$\%$ & ~ & 50$\%$ & 12 & 1202 & 0.071 & 94.9023$\%$\\\hline	
				\multirow{5}{*}{w6a} & 1$\%$ & 24 & 2712 & 0.165 & 98.3624$\%$ & \multirow{5}{*}{w7a} & 1$\%$ & 23 & 3138 & 0.193 & 98.6750$\%$\\
				~ & 5$\%$ & 7 & 959 & 0.073 & 96.3828$\%$ & ~ & 5$\%$ & 7 & 1144 & 0.096 & 97.0707$\%$\\
				~ & 10$\%$ & 8 & 1057 & 0.077 & 96.3980$\%$ & ~ & 10$\%$ & 9 & 1404 & 0.104 & 97.0627$\%$\\
				~ & 15$\%$ & 9 & 1199 & 0.085 & 96.3714$\%$ & ~ & 15$\%$ & 11 & 1653 & 0.116 & 97.0348$\%$\\
				~ & 25$\%$ & 9 & 1154 & 0.085 & 94.8060$\%$ & ~ & 25$\%$ & 9 & 1422 & 0.111 & 95.6579$\%$\\
				~ & 50$\%$ & 12 & 1477 & 0.106 & 93.9169$\%$ & ~ & 50$\%$ & 12 & 1820 & 0.135 & 94.8597$\%$\\\hline	
				\multirow{5}{*}{w8a} & 1$\%$ & 25 & 4586 & 0.350 & 98.6222$\%$ & \multirow{5}{*}{australian} & 1$\%$ & 20 & 1071 & 0.020 & 83.4025$\%$\\
				~ & 5$\%$ & 6 & 1346 & 0.160 & 97.0236$\%$ & ~ & 5$\%$ & 25 & 1305 & 0.023 & 83.6100$\%$\\
				~ & 10$\%$ & 9 & 1892 & 0.187 & 97.0303$\%$ & ~ & 10$\%$ & 23 & 1257 & 0.022 & 83.6100$\%$\\
				~ & 15$\%$ & 10 & 2084 & 0.201 & 96.9701$\%$ & ~ & 15$\%$ & 21 & 984 & 0.017 & 82.7801$\%$\\
				~ & 25$\%$ & 9 & 1918 & 0.192 & 95.6725$\%$ & ~ & 25$\%$ & 16 & 923 & 0.017 & 80.7054$\%$\\
				~ & 50$\%$ & 12 & 2467 & 0.227 & 95.0505$\%$ & ~ & 50$\%$ & 11 & 648 & 0.012 & 85.2697$\%$\\
			 \hline	
			\end{tabular}
		\end{center}
	}
\end{table}

\begin{table}[tbhp]
{\footnotesize
	\caption{{Numerical Results for \texttt{MPM} with Different Sparse Ratio (II).}} \label{tab:test_algo2}
	\begin{center}
			\begin{tabular}{p{20mm}p{2.5mm}p{2mm}p{1.5mm}p{5mm}p{14mm}|p{15mm}p{2.5mm}p{2mm}p{1.5mm}p{5mm}p{12mm}} \hline
				data set & $SR$ & k & $cg$ & t(s) & accuracy & data set & $SR$ & k & $cg$ & t(s) & accuracy  \\ \hline
					\multirow{5}{*}{breast-cancer} & 1$\%$ & 27 & 0 & 0.031 & 99.5122$\%$ & \multirow{5}{*}{cod-rna} & 1$\%$ & 23 & 0 & 0.074 & 60.7133$\%$\\
					~ & 5$\%$ & 25 & 0 & 0.001 & 99.5122$\%$ & ~ & 5$\%$ & 24 & 0 & 0.079 & 64.5037$\%$\\
					~ & 10$\%$ & 18 & 0 & 0.001 & 100.0000$\%$ & ~ & 10$\%$ & 22 & 0 & 0.073 & 64.7444$\%$\\
					~ & 15$\%$ & 15 & 0 & 0.001 & 99.5122$\%$ & ~ & 15$\%$ & 20 & 0 & 0.068 & 64.4701$\%$\\
					~ & 25$\%$ & 13 & 0 & 0.001 & 98.5366$\%$ & ~ & 25$\%$ & 16 & 0 & 0.060 & 62.0570$\%$\\
					~ & 50$\%$ & 10 & 0 & 0.001 & 99.0244$\%$ & ~ & 50$\%$ & 11 & 0 & 0.051 & 46.6379$\%$\\\hline	
					\multirow{5}{*}{colon-cancer} & 1$\%$ & 23 & 0 & 0.198 & 73.6842$\%$ & \multirow{5}{*}{diabetes} & 1$\%$ & 4 & 0 & 0.002 & 78.3550$\%$\\
					~ & 5$\%$ & 3 & 0 & 0.183 & 63.1579$\%$ & ~ & 5$\%$ & 7 & 0 & 0.001 & 78.7879$\%$\\
					~ & 10$\%$ & 3 & 0 & 0.183 & 63.1579$\%$ & ~ & 10$\%$ & 11 & 0 & 0.001 & 77.4892$\%$\\
					~ & 15$\%$ & 3 & 0 & 0.181 & 57.8947$\%$ & ~ & 15$\%$ & 16 & 0 & 0.001 & 77.0563$\%$\\
					~ & 25$\%$ & 3 & 0 & 0.180 & 57.8947$\%$ & ~ & 25$\%$ & 15 & 0 & 0.001 & 77.0563$\%$\\
					~ & 50$\%$ & 3 & 0 & 0.182 & 42.1053$\%$ & ~ & 50$\%$ & 14 & 0 & 0.001 & 67.5325$\%$\\\hline
				\multirow{5}{*}{duke breast-cancer} & 1$\%$ & 25 & 0 & 2.273 & 100.0000$\%$ & \multirow{5}{*}{fourclass} & 1$\%$ & 5 & 0 & 0.001 & 75.2896$\%$\\
				~ & 5$\%$ & 3 & 0 & 2.274 & 100.0000$\%$ & ~ & 5$\%$ & 11 & 0 & 0.001 & 76.4479$\%$\\
				~ & 10$\%$ & 3 & 0 & 2.348 & 100.0000$\%$ & ~ & 10$\%$ & 18 & 0 & 0.001 & 77.2201$\%$\\
				~ & 15$\%$ & 3 & 0 & 2.456 & 100.0000$\%$ & ~ & 15$\%$ & 22 & 0 & 0.005 & 76.0618$\%$\\
				~ & 25$\%$ & 3 & 0 & 2.473 & 100.0000$\%$ & ~ & 25$\%$ & 24 & 0 & 0.002 & 74.9035$\%$\\
				~ & 50$\%$ & 3 & 0 & 2.447 & 50.0000$\%$ & ~ & 50$\%$ & 11 & 0 & 0.001 & 72.2008$\%$\\\hline	
				\multirow{5}{*}{german.numer} & 1$\%$ & 5 & 0 & 0.002 & 80.0000$\%$ & \multirow{5}{*}{heart} & 1$\%$ & 6 & 0 & 0.001 & 87.6543$\%$\\
				~ & 5$\%$ & 10 & 0 & 0.004 & 80.0000$\%$ & ~ & 5$\%$ & 14 & 0 & 0.001 & 86.4198$\%$\\
				~ & 10$\%$ & 16 & 0 & 0.004 & 79.0000$\%$ & ~ & 10$\%$ & 17 & 0 & 0.003 & 86.4198$\%$\\
				~ & 15$\%$ & 21 & 0 & 0.004 & 76.0000$\%$ & ~ & 15$\%$ & 18 & 0 & 0.002 & 82.7160$\%$\\
				~ & 25$\%$ & 19 & 0 & 0.004 & 77.0000$\%$ & ~ & 25$\%$ & 14 & 0 & 0.001 & 87.6543$\%$\\
				~ & 50$\%$ & 13 & 0 & 0.004 & 73.0000$\%$ & ~ & 50$\%$ & 12 & 0 & 0.001 & 86.4198$\%$\\\hline	
				\multirow{5}{*}{ijcnn1} & 1$\%$ & 31 & 0 & 0.171 & 91.9565$\%$ & \multirow{5}{*}{ionosphere} & 1$\%$ & 23 & 0 & 0.003 & 95.2830$\%$\\
				~ & 5$\%$ & 30 & 0 & 0.184 & 91.7515$\%$ & ~ & 5$\%$ & 21 & 0 & 0.003 & 95.2830$\%$\\
				~ & 10$\%$ & 10 & 0 & 0.100 & 90.4996$\%$ & ~ & 10$\%$ & 20 & 0 & 0.002 & 98.1132$\%$\\
				~ & 15$\%$ & 11 & 0 & 0.106 & 90.5083$\%$ & ~ & 15$\%$ & 17 & 0 & 0.002 & 97.1698$\%$\\
				~ & 25$\%$ & 11 & 0 & 0.106 & 90.5137$\%$ & ~ & 25$\%$ & 13 & 0 & 0.002 & 97.1698$\%$\\
				~ & 50$\%$ & 14 & 0 & 0.114 & 87.9761$\%$ & ~ & 50$\%$ & 10 & 0 & 0.002 & 98.1132$\%$\\\hline	
				\multirow{5}{*}{leukemia} & 1$\%$ & 36 & 0 & 2.265 & 64.7059$\%$ & \multirow{5}{*}{liver-disorders} & 1$\%$ & 4 & 0 & 0.000 & 50.0000$\%$\\
				~ & 5$\%$ & 3 & 0 & 2.236 & 64.7059$\%$ & ~ & 5$\%$ & 6 & 0 & 0.000 & 50.0000$\%$\\
				~ & 10$\%$ & 3 & 0 & 2.255 & 61.7647$\%$ & ~ & 10$\%$ & 8 & 0 & 0.000 & 50.0000$\%$\\
				~ & 15$\%$ & 3 & 0 & 2.240 & 76.4706$\%$ & ~ & 15$\%$ & 9 & 0 & 0.000 & 50.0000$\%$\\
				~ & 25$\%$ & 3 & 0 & 2.242 & 82.3529$\%$ & ~ & 25$\%$ & 11 & 0 & 0.001 & 50.0000$\%$\\
				~ & 50$\%$ & 2 & 0 & 2.218 & 79.4118$\%$ & ~ & 50$\%$ & 16 & 0 & 0.001 & 50.0000$\%$\\\hline	
				\multirow{5}{*}{mushrooms} & 1$\%$ & 25 & 767 & 0.013 & 100.0000$\%$ & \multirow{5}{*}{phishing} & 1$\%$ & 23 & 0 & 0.054 & 91.2270$\%$\\
				~ & 5$\%$ & 18 & 569 & 0.010 & 100.0000$\%$ & ~ & 5$\%$ & 23 & 0 & 0.051 & 92.0711$\%$\\
				~ & 10$\%$ & 14 & 494 & 0.008 & 100.0000$\%$ & ~ & 10$\%$ & 13 & 0 & 0.043 & 89.3579$\%$\\
				~ & 15$\%$ & 15 & 535 & 0.010 & 100.0000$\%$ & ~ & 15$\%$ & 11 & 0 & 0.041 & 88.9358$\%$\\
				~ & 25$\%$ & 14 & 562 & 0.010 & 100.0000$\%$ & ~ & 25$\%$ & 11 & 0 & 0.045 & 89.4784$\%$\\
				~ & 50$\%$ & 12 & 429 & 0.013 & 99.3590$\%$ & ~ & 50$\%$ & 12 & 0 & 0.046 & 90.4432$\%$\\\hline	
				\multirow{5}{*}{splice} & 1$\%$ & 12 & 0 & 0.009 & 64.0920$\%$ & \multirow{5}{*}{spnar} & 1$\%$ & 22 & 0 & 0.004 & 17.4603$\%$\\
				~ & 5$\%$ & 24 & 0 & 0.011 & 62.2989$\%$ & ~ & 5$\%$ & 21 & 0 & 0.002 & 23.8095$\%$\\
				~ & 10$\%$ & 25 & 0 & 0.010 & 68.8276$\%$ & ~ & 10$\%$ & 18 & 0 & 0.002 & 31.7460$\%$\\
				~ & 15$\%$ & 22 & 0 & 0.011 & 71.6322$\%$ & ~ & 15$\%$ & 20 & 0 & 0.002 & 25.3968$\%$\\
				~ & 25$\%$ & 20 & 0 & 0.011 & 64.5057$\%$ & ~ & 25$\%$ & 13 & 0 & 0.002 & 31.7460$\%$\\
				~ & 50$\%$ & 13 & 0 & 0.009 & 50.2069$\%$ & ~ & 50$\%$ & 14 & 0 & 0.002 & 33.3333$\%$\\\hline	
				\multirow{5}{*}{svmguide1} & 1$\%$ & 29 & 0 & 0.007 & 95.2750$\%$ & \multirow{5}{*}{svmguide3} & 1$\%$ & 8 & 0 & 0.003 & 43.9024$\%$\\
				~ & 5$\%$ & 28 & 0 & 0.007 & 93.9250$\%$ & ~ & 5$\%$ & 13 & 0 & 0.005 & 48.7805$\%$\\
				~ & 10$\%$ & 26 & 0 & 0.006 & 92.1750$\%$ & ~ & 10$\%$ & 13 & 0 & 0.004 & 43.9024$\%$\\
				~ & 15$\%$ & 24 & 0 & 0.006 & 90.3000$\%$ & ~ & 15$\%$ & 13 & 0 & 0.004 & 34.1463$\%$\\
				~ & 25$\%$ & 21 & 0 & 0.006 & 88.7000$\%$ & ~ & 25$\%$ & 10 & 0 & 0.003 & 2.4390$\%$\\
				~ & 50$\%$ & 11 & 0 & 0.003 & 74.9000$\%$ & ~ & 50$\%$ & 13 & 0 & 0.004 & 60.9756$\%$\\ \hline
			\end{tabular}
		\end{center}
	}
\end{table}

\subsection{Numerical Comparisions with LIBLINEAR}

{In this part, we compare our algorithm with some solvers in \texttt{LIBLINEAR} \cite{fan2008liblinear}. \texttt{LIBLINEAR} is an open source library for large-scale linear classification. We choose the dual coordinate descent method (\texttt{DCD}) in \cite{hsieh2008dual} and the trust-region Newton (\texttt{TRN}) in \cite{lin2007trust}. \texttt{DCD} is designed for the dual problem of \texttt{L1-SVM} and \texttt{L2-SVM}. \texttt{TRN} is used {for the primal from of} \texttt{L2-SVM}. 

For our \texttt{MPM}, the parameters are consistent {as} in Section 5.1. We set the $SR$ as $10\%$. For \texttt{DCD} and \texttt{TRN}, all parameters {are set as} default. 
We report the cputime and accuracy for each algorithm in Table \ref{tab:compare}.

\begin{table}[tbhp]
	{\footnotesize
		\caption{ The Comparison Results. A1: \texttt{MPM}; A2: \texttt{DCD} for dual \texttt{L2-SVM};
			A3: \texttt{TRN} for primal \texttt{L2-SVM}; A4: \texttt{DCD} for dual \texttt{L1-SVM}. 
		} \label{tab:compare} 
		\begin{center}
		
				\begin{tabular}{ccc} \hline
					data set & t(s) (A1$|$A2$|$A3$|$A4) & accuracy (A1$|$A2$|$A3$|$A4) \\ \hline
					a1a & 0.029$|$0.040$|$0.002$|$0.004 & 83.50$|$83.84$|$83.85$|$83.81 \\
					a2a & 0.033$|$0.012$|$0.003$|$0.006 & 84.19$|$84.00$|$84.04$|$84.28 \\
					a3a & 0.042$|$0.018$|$0.003$|$0.008 & 83.71$|$84.26$|$84.32$|$84.29 \\
					a4a & 0.059$|$0.028$|$0.005$|$0.016 & 84.33$|$84.41$|$84.40$|$84.49 \\
					a5a & 0.066$|$0.039$|$0.008$|$0.020 & 84.37$|$84.52$|$84.51$|$84.38 \\
					a6a & 0.094$|$0.076$|$0.013$|$0.037 & 84.54$|$84.80$|$84.79$|$84.66 \\
					a7a & 0.119$|$0.115$|$0.021$|$0.060 & 84.76$|$85.04$|$85.01$|$84.85 \\
					a8a & 0.148$|$0.187$|$0.031$|$0.119 & 85.27$|$85.44$|$85.38$|$85.10 \\
					a9a & 0.228$|$0.315$|$0.041$|$0.142 & 84.73$|$84.99$|$84.99$|$85.00 \\
					w1a & 0.024$|$0.003$|$0.002$|$0.002 & 95.17$|$95.98$|$95.94$|$96.00 \\
					w2a & 0.105$|$0.214$|$0.032$|$0.062 & 95.99$|$97.95$|$97.96$|$98.00 \\
					w3a & 0.023$|$0.003$|$0.001$|$0.001 & {\bf97.51}$|$94.78$|$94.66$|$93.99 \\
					w4a & 0.037$|$0.016$|$0.007$|$0.006 & 94.77$|$96.77$|$96.76$|$96.70 \\
					w5a & 0.057$|$0.033$|$0.010$|$0.014 & 96.97$|$97.96$|$97.95$|$97.95 \\
					w6a & 0.075$|$0.096$|$0.019$|$0.023 & 96.40$|$97.88$|$97.87$|$97.83 \\
					w7a & 0.106$|$0.129$|$0.026$|$0.046 & 97.06$|$98.19$|$98.17$|$98.18 \\
					w8a & 0.190$|$0.438$|$0.060$|$0.100 & 97.03$|$98.34$|$98.37$|$98.39 \\
					australian & 0.023$|$0.005$|$0.002$|$0.003 & {\bf83.61}$|$83.40$|$83.20$|$83.20 \\
					breast-cancer & 0.003$|$0.001$|$0.000$|$0.000 & {\bf100.00}$|$99.51$|$99.51$|$99.51 \\
					{\bf cod-rna} & 0.073$|$2.867$|$0.027$|$2.929 & {\bf64.74}$|$47.84$|$58.21$|$48.03 \\
					\underline{colon-cancer} & 0.186$|$0.003$|$0.003$|$0.003 & 63.16$|$68.42$|$68.42$|$68.42 \\
					diabetes & 0.003$|$0.001$|$0.000$|$0.001 & 77.49$|$79.22$|$79.22$|$77.49 \\
					\underline{duke breast-cancer} & 2.250$|$0.013$|$0.015$|$0.014 & {\bf100.00}$|$100.00$|$100.00$|$100.00 \\
					{\bf fourclass} & 0.003$|$0.001$|$0.000$|$0.001 & {\bf77.22}$|$65.25$|$65.64$|$69.11 \\
					german.numer & 0.005$|$0.008$|$0.001$|$0.003 & 79.00$|$79.33$|$79.67$|$79.00 \\
					heart & 0.001$|$0.001$|$0.000$|$0.000 & {\bf86.42}$|$83.95$|$83.95$|$83.95 \\
					ijcnn1 & 0.100$|$0.123$|$0.062$|$0.082 & 90.50$|$91.79$|$91.78$|$92.11 \\
					ionosphere & 0.005$|$0.003$|$0.001$|$0.001 & 98.11$|$98.11$|$98.11$|$99.06 \\
					\underline{leukemia} & 2.272$|$0.011$|$0.014$|$0.011 & 61.76$|$79.41$|$79.41$|$79.41 \\
					liver-disorders & 0.001$|$0.000$|$0.000$|$0.000 & {\bf50.00}$|$50.00$|$50.00$|$50.00 \\
					mushrooms & 0.011$|$0.001$|$0.001$|$0.001 & {\bf100.00}$|$100.00$|$100.00$|$100.00 \\
					phishing & 0.044$|$0.015$|$0.015$|$0.012 & 89.36$|$91.59$|$91.59$|$91.62 \\
					{\bf splice} & 0.011$|$0.055$|$0.003$|$0.014 & {\bf68.83}$|$52.05$|$52.09$|$57.56 \\
					{\bf spnar} & 0.004$|$0.002$|$0.001$|$0.002 & {\bf31.75}$|$11.11$|$11.11$|$11.11 \\
					{\bf svmguide1} & 0.006$|$0.103$|$0.001$|$0.111 & {\bf92.17}$|$78.55$|$78.92$|$62.65 \\
					{\bf svmguide3} & 0.003$|$0.009$|$0.002$|$0.003 & {\bf43.90}$|$24.39$|$24.39$|$9.76\\\hline				
				\end{tabular}
			\end{center}
		}
	\end{table}

The results obtained in Table \ref{tab:compare} are {analyzed} as follows.
The {four} methods {provide} good accuracy {for most the data sets}. {For} most of the data sets, more than 70$\%$, {or} even 90$\%$ {of the test data} {can be correctly predicted with right labels.} But for some very difficult data sets, our algorithm has a very good {result. That is, there is significant improvement in accuracy for solving} cod-rna, fourclass, splice, spnar, svmguide1 and svmguide3 {(indicated by bold font in Table \ref{tab:compare}).} However, when the feature dimension $m$ of {the} individual data set is much larger than the sample number $n$  ({ colon-cancer, duke breast-cancer and leukemia} {indicated by underline}), \texttt{MPM} takes a bit more time. {This can be explained by the fact that} 
{the resulting subproblem is a $m$ by $m$ linear system, which takes relatively more time to solve than other data sets.}

\section{Conclusions}

In our paper, we proposed our \texttt{SCSVM} model with 0-norm constraint, which is different from \texttt{L1-SVM} and \texttt{L2-SVM}. We penalized the sparse constraint to the objective function and used the majorizaiton method to solve it approximately. It should be noted that our method itself is not limited to our SVM model. We proposed Theorem \ref{theorem1.1} {to deal with more general cases}, which shows that as long as {conditions (a) and (b) are satisfied,} {we can apply the majorization penalty method.}

{The subproblem in the majorization penalty method admits a special form, which is a strongly convex quadratic problem. Therefore, it is equivalent to solving a linear system.} We proved the nonsingularity of the Jacobian matrix, and then we used {conjugate gradient method} to solve {the linear} equation. In the numerical experiments, we compared with two methods \texttt{DCD} and \texttt{TRN} for models \texttt{L1-SVM} and \texttt{L2-SVM} respectively in LIBLINEAR. {The results verified the efficiency of the proposed model and our algorithm.}

\bibliographystyle{plain}
\bibliography{references_ssc}

\end{document}